\documentclass[conference,10pt]{IEEEtran}

\usepackage{cite}
\usepackage{amsmath,amssymb,amsfonts}
\usepackage{algorithm,algcompatible,algpseudocode}
\usepackage{graphicx}
\usepackage{textcomp}
\usepackage{bm}
\usepackage{float}
\usepackage{longtable,supertabular}
\usepackage{booktabs,chemformula}
\usepackage{url}
\usepackage{microtype}
\usepackage{cite}
\usepackage{balance}
\usepackage[hidelinks]{hyperref}

\renewcommand{\eqref}[1]{(\ref{#1})}

\newcommand{\setI}{\ensuremath{\mathcal{I}}}

\newcommand{\bma}{\ensuremath{\mathbf{a}}}
\newcommand{\bmb}{\ensuremath{\mathbf{b}}}

\newcommand{\bmx}{\ensuremath{\mathbf{x}}}

\newcommand{\bA}{\ensuremath{\mathbf{A}}}

\newcommand{\bI}{\ensuremath{\mathbf{I}}}

\newcommand{\bY}{\ensuremath{\mathbf{Y}}}

\newcommand{\aln}[1]{\begin{align}#1\end{align}}

\newcommand{\matb}{\left( \begin{matrix*}[r] }
\newcommand{\mate}{\end{matrix*}\right)}

\begin{document}

\title{PhasePack: A Phase Retrieval Library\\
{\large \url{https://github.com/tomgoldstein/phasepack-matlab}}
\thanks{PhasePack was made possible by the NSF REU program (CCF-1535902).  The work of T.~Goldstein was supported in part by the US National Science Foundation (NSF) under grant CCF-1535902, the US Office of Naval Research under grant N00014-17-1-2078, and by the Sloan Foundation. The work of C. Studer was supported in part by Xilinx, Inc. and by the US NSF under grants ECCS-1408006, CCF-1535897,  CNS-1717559, and CAREER CCF-1652065.}
}

\author{\IEEEauthorblockN{Rohan Chandra\IEEEauthorrefmark{1},
Ziyuan Zhong\IEEEauthorrefmark{2}, 
Justin Hontz\IEEEauthorrefmark{1}, 
Val McCulloch\IEEEauthorrefmark{3}, 
Christoph Studer\IEEEauthorrefmark{4}   and
Tom Goldstein\IEEEauthorrefmark{1}\\
 \IEEEauthorrefmark{1}University of Maryland,
\IEEEauthorrefmark{2}Columbia University,
\IEEEauthorrefmark{3}Smith College,
\IEEEauthorrefmark{4}Cornell University}}
\maketitle

%
%
%
%
%
%
%
%
%
%

\maketitle

\begin{abstract}
Phase retrieval deals with the estimation of complex-valued signals solely from the magnitudes of linear measurements. While there has been a recent explosion in the development of phase retrieval algorithms, the lack of a common interface has made it difficult to compare new methods against the state-of-the-art. The purpose of PhasePack is to create a common software interface for a wide range of phase retrieval algorithms and to provide a common testbed using both synthetic data and empirical imaging datasets. 
PhasePack is able to benchmark a large number of recent phase retrieval methods against one another to generate comparisons using a range of different performance metrics. The software package handles single method testing as well as multiple method comparisons.

The algorithm implementations in PhasePack differ slightly from their original descriptions in the literature in order to achieve faster speed and improved robustness. In particular,  PhasePack uses adaptive stepsizes, line-search methods, and fast eigensolvers to speed up and automate convergence.
 \end{abstract}
\smallskip

\section{Introduction}


Classical linear inverse problems recover an unknown signal $\bmx\in \mathbb{R}^n$ from $m$ linear measurements of the form $\bA \bmx=\bmb,$ where $\bA\in \mathbb{C}^{m\times n}$ is a (known) linear measurement operator and $\bmb\in \mathbb{C}^m$ is a (possibly noisy) vector of measurements.   In contrast, phase retrieval \cite{fienup1982phase,jaganathan2015phase,shechtman2015phase} is the recovery of signals from the {\em magnitudes} of linear measurements $|\bA \bmx|, $ where $|\cdot|$ denotes the entry-wise absolute value of the (possibly complex-valued) measurements.  
%
Phase retrieval is prevalent in the physical sciences where detectors can measure only the magnitude (brightness or intensity) of an electromagnetic wave, but not its phase. Such modalities include optical imaging \cite{millane1990phase,bian2015fourier}, quantum state tomography, electron microscropy \cite{misell1973method}, X-ray crystallography, astronomy \cite{fienup1982phase}, and X-ray diffraction imaging \cite{shechtman2015phase}.
%
%
%

Phase Retrieval is an example of a non-convex quadratic program with quadratic constraints.  In the real-valued case, it is a combinatorial problem of determining the missing signs of $\bA\bmx$, which is known to be NP-hard \cite{pardalos1991quadratic}.
Despite this observation, recent years have seen the development of new algorithms that solve phase retrieval problems effectively.  Unfortunately, because of the lack of publicly available real-world data, the lack of a common software interface for different algorithms, and a knowledge gap between practitioners and theoreticians, only little work has been devoted to compare and evaluate newer phase retrieval methods. 

\section{Introducing PhasePack} 
PhasePack is a software package that contains implementations of many different phase retrieval methods, and tools for easily applying them to real and synthetic datasets.  The purpose of PhasePack is to create a common interface for a wide range of phase retrieval schemes, and to provide a common testbed using both synthetic data and empirical imaging datasets \cite{metzler2017coherent}. 
PhasePack can also benchmark different algorithms against one another, and generate performance comparisons with varying numbers of measurements, signal-to-noise ratio, iterations, and running time. The package handles single method testing as well as multiple method comparisons.

The methods in PhasePack differ in numerous ways from their original descriptions in the literature in order to achieve improved robustness to different measurement models and to enable faster convergence.  After reviewing the kinds of methods available in PhasePack in Section \ref{methods}, we will discuss practical implementation consideration in Section \ref{practical}.

\section{Methods available in PhasePack}\label{methods}
The methods implemented in PhasePack can be divided into three categories.  First, alternating minimization methods work by iteratively estimating the missing phases of the measurements, and then solving the linear system $\bA\bmx = \tilde \bmb$ where $\tilde \bmb$ is formed by applying phase estimates to $\bmb.$ This class of algorithms includes the classical Gerchberg-Saxton \cite{gerchberg1972practical} and Fienup \cite{fienup1982phase} methods, in addition to more recent methods proposed in~\cite{netrapalli2013phase,marchesini2016alternating,chen2015phase}.
%

Second, many approaches attack the non-convex formulation directly using least-squares formulations \cite{candes2015wirtinger,chen2015solving,wang2016solving,yuan2017phase,candes2015phase,wang2017solving, zeng2017coordinate,chen2015solving,li2016gradient,wei2015solving}. These methods use gradient descent to  minimize an objective of the form
\aln{ \label{ls}
  f(\bmx)=  \| \, |\hat\bA\bmx|^p - \hat \bmb^p \, \|_2
}
where $\hat \bA$ and $\hat \bmb$ contain either a subset of the measurements or a re-weighting of the system $\bA\bmx=\bmb,$ and $p\in \{1,2\}$ is an integer exponent.  Because of the non-convexity of the formulation \eqref{ls}, these methods require careful initialization to avoid local minimizers.

The third class of methods convert phase retrieval into a convex problem.  This includes the ``lifting methods'' PhaseLift~\cite{candes2013phaselift} and PhaseCut \cite{waldspurger2015phase}, which square the dimensionality of the problem.  PhasePack also implements several low-rank approximate solvers for large instances of PhaseLift, including gauge duality methods \cite{friedlander2014gauge,aravkin2017foundations} and sketching methods~\cite{yurtsever2017sketchy}.  Finally, PhasePack contains the newer non-lifting relaxations PhaseMax~\cite{bahmani2017phase,goldstein2017icml,goldstein2016phasemax} and PhaseLamp \cite{dhifallah2017phase}.

\section{Building practical implementations} \label{practical}
\begin{figure}[tp]
\centering
\includegraphics[width=9cm, trim=0cm 7cm 1cm 7.0cm, clip]{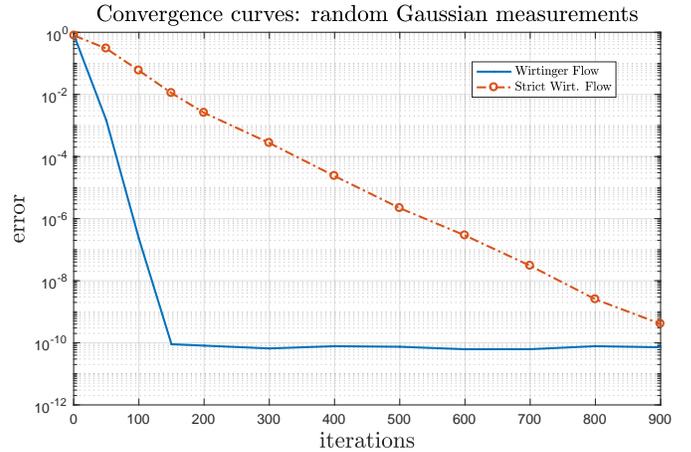}
\caption{ Reconstruction of a signal using two variants of Wirtinger flow. (blue) Optimized Wirtinger flow with adaptive stepsize and backtracking line search. (red) A ``strict'' implementation of Wirtinger flow using the non-adaptive stepsize rule described in the original paper. }
\label{fig:strict}
\end{figure}

Many recent phase retrieval methods have been developed with the goal of proving rigorous guarantees for random Gaussian measurement matrices.  As a result, the implementations described in the literature are not optimized for performance, and may become unstable when used on non-Gaussian measurements.  This is particularly an issue for variations of Wirtinger flow, which require the choice of stepsize parameters.  This is illustrated in Figure \ref{fig:strict}, which compares a ``strict'' implementation of Wirtinger flow using the stepsize rules described  in \cite{candes2015phase} to PhasePack's implementation with adaptive stepsize.

The implementations in PhasePack depart from the literature in several ways that make them more robust and efficient.  We discuss several such issues below.

\subsection{Adaptive stepsizes and backtracking}

 Wirtinger Flow and its variants rely on gradient descent methods, and are sensitive to the choice of stepsize parameters.  Existing stepsize rules presented in the literature are designed for Gaussian measurement models, and may be unstable or slow when used on generic matrices.
PhasePack solves least-squares formulations of the phase retrieval problem using the general gradient descent solver FASTA \cite{FASTA:2014}, which supports adaptive stepsizes, automated stopping conditions, and conjugate-gradient acceleration.  

This gradient solver automates optimization in several ways.  For stepsize selection, PhasePack uses the Barzilai-Borwein adaptive method \cite{BB88}.  
To guarantee stability, the method uses a backtracking line-search similar to the classical Armijo line search \cite{BMDIC95}.  Classical line searches enforce that the objective decreases monotonically.  For non-convex problems where the local curvature changes rapidly and local minima are prevalent, monotonic searches may result in excessive backtracking, and may be more prone to getting stuck in bad local minimizers.    For this reason, several authors have proposed non-monotonic line searches \cite{GLL86,ZH04}. Rather than requiring $f(\bmx^{(k+1)})<f(\bmx^{(k)})$ on the $k$th iteration, these methods require the weaker condition $f(\bmx^{(k+1)}) < \max_{k-w<j<k}\{f(\bmx^{(j)})\}$ for some ``window'' parameter $w.$    The particular adaptive strategy and line search rule used by PhasePack is described in detail in \cite{FASTA:2014}. 

Finally, the solver in PhasePack is capable of using L-BFGS acceleration (which was studied for phase retrieval in \cite{li2016gradient}) and non-linear conjugate gradient methods to speed up convergence. However, these options need to be tuned manually.

\begin{figure*}[tp]
\centering
\includegraphics[width=9cm, trim=0cm 8cm 0cm 8.0cm, clip]{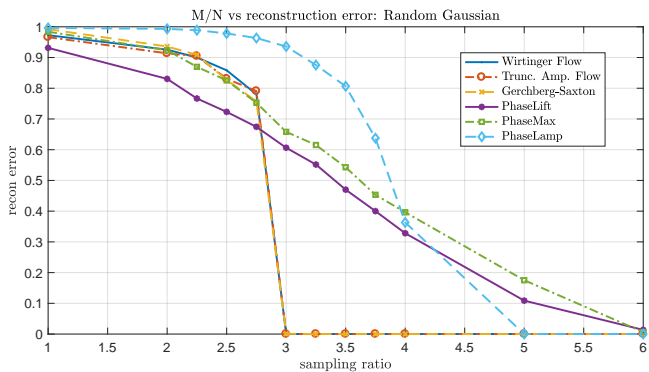}~
\includegraphics[width=9cm, trim=0cm 8cm 0cm 8.0cm, clip]{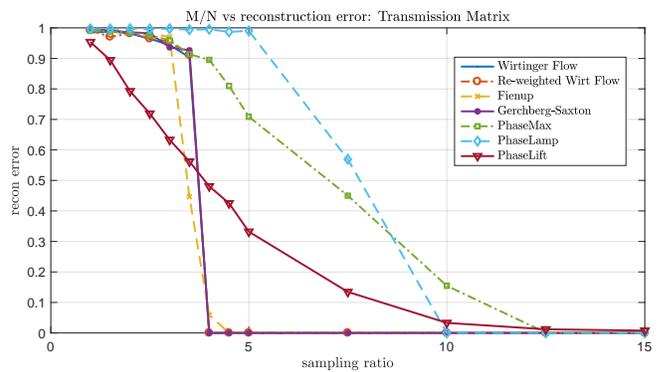}
\caption{ Reconstruction error as a function of sampling ratio (number of measurements over number of unknowns) for a few algorithms. (left)
Signal of length 100 using synthetic Gaussian measurements. (right) Signal of length 256 using an empirical transmission matrix measurement operator.}
\label{fig:random}
\end{figure*}

\subsection{Practical methods for spectral initialization}
\label{init}
Non-convex methods for phase retrieval are prone to getting stuck in local minimizers if they are not initialized properly.  Phasepack includes two main classes of initializers: (i) the spectral initializer and its variants and, (ii) the orthogonality promoting initializer.  The initializers described in the literature are optimized for random Gaussian data, and may not perform well on non-Gaussian real-world data.   To be robust to different measurement models, PhasePack initializers depart from the literature in three main ways:  we use a Krylov subspace-based eigensolver, we re-scale data before pre-processing, and we re-scale the initialization vector after computation.
  We discuss these important differences below after briefly describing how initialization methods work.

Spectral initialization methods begin by forming the matrix
\aln{\label{spectral}
 \bY = \dfrac{1}{m} \sum_{i=1}^m  T(b_i^2)  \bma_i  \bma_i^H 
 }
 where $\bma_i$ is the $i$th row of $\bA$, and $T(\cdot)$ is a ``pre-processing'' function, which is simply the identity in the original spectral method \cite{netrapalli2013phase}.  These methods are motivated by the observation that if $\bma$ is a random Gaussian vector and $b$ is its corresponding measurement, then we have
\aln{
\mathbb{E}_\bma [b^2  \bma  \bma^H] = 2 \bmx  {\bmx^H} + \| \bmx\|^2 \bI_n.
\label{expect}
}
In fact, the leading eigenvector of \eqref{expect} is the unknown signal~$\bmx.$ Spectral methods work by approximating the expectation~\eqref{expect} with the empirical matrix \eqref{spectral}, computing the leading eigenvector, and then scaling it appropriately.  A variety of methods have been proposed with different pre-processing functions, including the identity \cite{netrapalli2013phase,candes2015phase}, a ``truncation'' operator that is zero for small magnitude measurements \cite{chen2015solving}, and a ``re-weighting'' operator that shrinks its argument \cite{yuan2017phase}.   A rigorous study of these methods appears in \cite{lu2017phase}.  The default initialization in PhaseMax uses the ``optimal'' preprocessing function proposed in \cite{mondelli2017fundamental}, which is 
  \aln{ \label{T}
  T(z) =\frac{z-1}{z+\sqrt{\delta}-1}
  }
  where $\delta=m/n$ is the number of measurements divided by the signal dimension.

%
%
%
%
%
%
%
%

Orthogonality promoting initializers \cite{wang2016solving} work by identifying measurement vectors that are {\em not} correlated with the signal $\bmx$, and finding an initialization vector that is orthogonal to those uncorrelated vectors.  
This is done by forming the matrix
\[ \bY = \dfrac{1}{|\setI|} \sum_{i \in  \setI} \dfrac{ \bma_i { \bma_i}^H}{\| \bma_i\|^2}\] 
where $\setI$ contains the indices of measurement vectors that produce the smallest measurements.  The initializer is then the eigenvector of $\bY$ with {\em smallest} eigenvalue.

PhasePack initializers depart from the literature in several ways that we now describe.   First the literature relies on the power method to compute leading eigenvectors.  This is limiting for several reasons.  The power method can be unstable if the matrix does not have a large spectral gap (which empirical matrices may have not), and it is incapable of finding the smallest eigenvector (which is needed for orthogonality promoting initializers).  For this reason, we compute eigenvectors using the iterated Arnoldi method with restart, which is the default algorithm behind MATLAB's ``eigs'' routine.  This method exhibits faster convergence than the power method and can reliably resolve both the leading and trailing eigenvectors. 

A more major departure from the literature is in how PhasePack scales data before pre-processing.  The pre-processing functions $T$ defined in the literature are optimized for random Gaussian measurement matrices and normalized signals, and initializers may be highly sub-optimal (or even ill-posed) when empirical data is used.  For this reason, data is re-scaled to have the same statistics as the Gaussian model before pre-processing.  Note that the measurement operator may be implemented as a black-box function, and so re-scaling needs to be done without access to the matrix entries or the underlying signal. 

 For example, the pre-processing operator \eqref{T} is optimized for Gaussian random measurements of variance $1/n$ and a signal of length $n$ \cite{mondelli2017fundamental}, and has unpredictable results for any other model/signal.  However, this particular measurement model produces measurements with $\mathbb{E}_\bma |\bma^H\bmx|^2=1.$   For this reason, instead of applying the pre-processing operator to $\bmb$, we apply the pre-processing operator to the re-scaled measurement vector $\frac{m\bmb}{\|\bmb\|},$ which has the same mean as the Gaussian model.

Finally, PhasePack concludes initialization with a least-squares step to determine the optimal length of the initial vector.  Most eigensolvers produce eigenvalues that have unit length by default, however, the length of the signal may effect the convergence of the phase retrieval solver.  For this reason, after finding the initializer $\hat \bmx$ we compute the least-squares solution to 
  $$\min_\alpha \|  \,  \alpha |\bA\hat \bmx| - \bmb   \, \|,$$
  and then replace the initializer $\hat \bmx$ with $\alpha \hat \bmx.$   The spectral initialization methods described above come with their own scaling routines, however PhasePack's least-squares scaling method is more practical for empirical data as it does not depend on the Gaussian measurement model.  

%
%
%

\begin{figure*}[th]
\centering
\begin{minipage}{3cm} \centering
\vspace{12pt}Original\vspace{14pt}
\includegraphics[width=3cm]{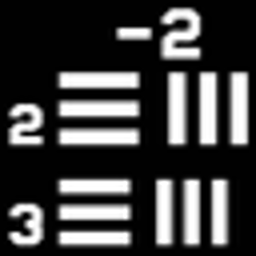}

\end{minipage}
\begin{minipage}{3cm} \centering
\vspace{12pt}Fienup \\
0.13019\vspace{3pt}
\includegraphics[width=3cm]{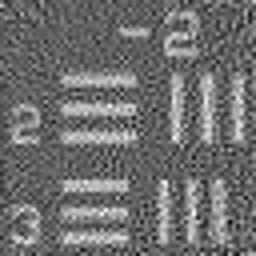}
\end{minipage}
\begin{minipage}{3cm} \centering
\vspace{12pt}Gerchberg-Saxton \\
0.13019\vspace{3pt}
\includegraphics[width=3cm]{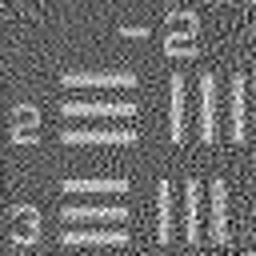}
\end{minipage}
\begin{minipage}{3cm} \centering
\vspace{12pt}Wirtinger Flow\\0.22900\vspace{3pt}
\includegraphics[width=3cm]{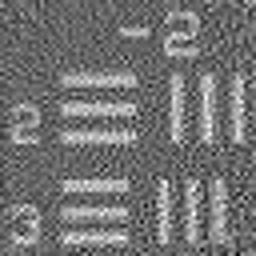}
\end{minipage}
\begin{minipage}{3cm} \centering
Truncated Wirtinger Flow\\0.18017\vspace{3pt}
\includegraphics[width=3cm]{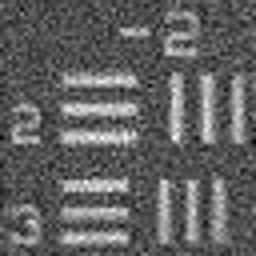}
\end{minipage} \\
\begin{minipage}{3cm} \centering
\vspace{3pt}Reweighted Wirtinger Flow\\0.14044\vspace{3pt}
\includegraphics[width=3cm]{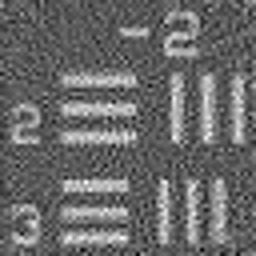}
\end{minipage} 
\begin{minipage}{3cm} \centering
\vspace{15pt}Amplitude Flow\\0.13017\vspace{3pt}
\includegraphics[width=3cm]{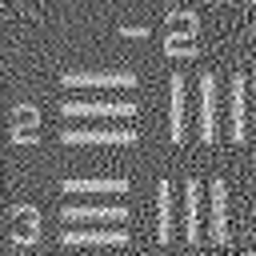}
\end{minipage}
\begin{minipage}{3cm} \centering
\vspace{3pt}Truncated Amplitude Flow\\0.43330\vspace{3pt}
\includegraphics[width=3cm]{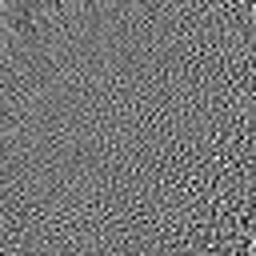}
\end{minipage}
\begin{minipage}{3cm} \centering
\vspace{3pt}Reweighted Amplitude Flow\\ 0.13261\vspace{3pt}
\includegraphics[width=3cm]{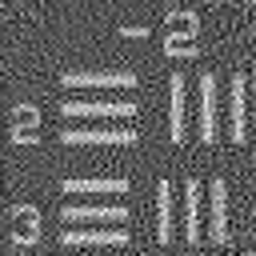}
\end{minipage}
\begin{minipage}{3cm} \centering
\vspace{15.25pt}Kaczmarz\\0.34500\vspace{3pt}
\includegraphics[width=3cm]{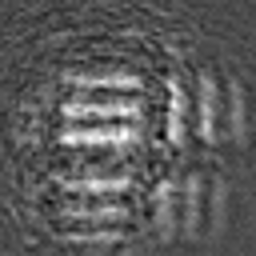}
\end{minipage}\\
\begin{minipage}{3cm} \centering
\vspace{17pt}Coordinate Descent\\0.21616\vspace{3pt}
\includegraphics[width=3cm]{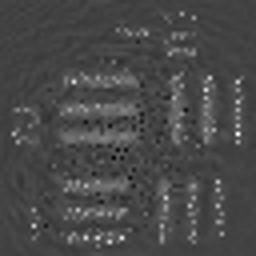}
\end{minipage}
\begin{minipage}{3cm} \centering
\vspace{17pt}PhaseLift\\0.35452\vspace{3pt}
\includegraphics[width=3cm]{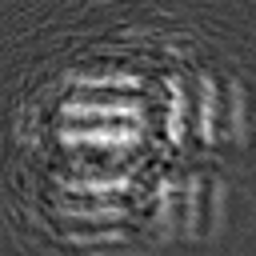}
\end{minipage}
\begin{minipage}{3cm} \centering
\vspace{17pt}PhaseMax\\0.23459\vspace{3pt}
\includegraphics[width=3cm]{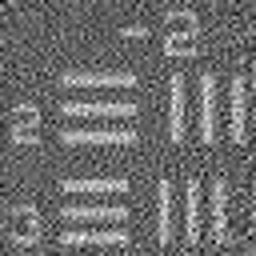}
\end{minipage}
\begin{minipage}{3cm} \centering
\vspace{17pt}PhaseLamp\\0.64980\vspace{3pt}
\includegraphics[width=3cm]{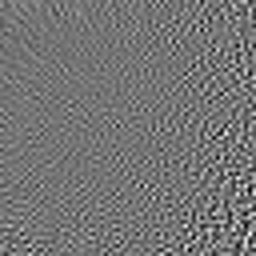}
\end{minipage}
\begin{minipage}{3cm} \centering
\vspace{17pt}SketchyCGM\\0.85526\vspace{3pt}
\includegraphics[width=3cm]{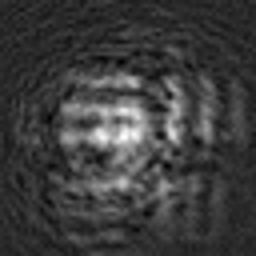}
\end{minipage}
\vspace{3mm}
\caption{ Reconstruction of a $64\times 64$ pixel image (modified air force target) using empirical phaseless measurements. Below each algorithm name is 
the relative measurement error, defined as $\||\bA\bmx|-\bmb\|/\|\bmb\|$, achieved by each method.} 
\label{fig:aft}
\end{figure*}

\section{Empirical datasets}
Public datasets for phase retrieval are scarce, and available datasets are often difficult to use because of lack of documentation or unusual data formats.  One of PhasePacks goals it to create a simple API for testing methods on empirical datasets.  
PhasePack provides routines for unpacking and pre-processing datasets, and prepares a simple measurement operator $\bA$ and measurement vector $\bmb$ that can be used for phase retrieval. 
 The datasets currently supported by PhasePack were obtained for the purpose of imaging an object through a diffusive medium, and are described in detail in \cite{metzler2017coherent}.   We provide a brief description of this imaging modality below.

The imaging modality described  \cite{metzler2017coherent} is for reconstructing an image from light that has passed through a diffusive medium (e.g., a light-scattering material like paper or eggshell).   Imaging though such a material is impossible using a conventional camera, but this can be done using phase retrieval methods.  A binary mask is created by placing a transparent LCD display in front of a light source.  The pixels in this mask can be made either opaque or clear.  An image pattern in loaded onto the mask, and the mask is illuminated from behind using a coherent (laser) light source. The light passes through the mask (creating an image), and then passes through a diffuser.  The diffused light then lands on a photodetector that measures the intensity of light (but not the phase) at many locations, producing measurements.  

  The light pattern (image) passing through the mask is the unknown signal $\bmx$.  The diffuser is modeled as a ``transmission matrix'' $\bA$ that describes how light passing through each pixel of the mask effects each of the detectors.  Under this model, the signal received at each detector is given by $\bA\bmx.$ Because standard photodetectors measure intensity and not phase, the recorded signal is $\bmb = |\bA\bmx|,$ and the image must be reconstructed by phase retrieval.
  
  This imaging modality is particularly useful for benchmarking phase retrieval methods because the ground-truth image (the pattern loaded onto the mask) is known, and can be used to evaluate the quality of reconstructions.  A sample mask and its reconstructions is shown in Figure \ref{fig:aft}.

\section{Experimental Results}
We demonstrate the capabilities of PhasePack using a range of experiments on real and synthetic data.
All of these examples were produced using the scripts included in the \path{benchmark} sub-folder of the PhasePack distribution.  Figure \ref{fig:random} shows the performance of several algorithms for recovering a random signal from random Gaussian measurements (left) and from measurements acquired using a real empirical transmission matrix (right).  Figure \ref{fig:aft} shows the reconstruction of an image from empirical measurements obtained using an optical device~\cite{metzler2017coherent}.  

In most of our experiments, we found that the classical Fienup \cite{fienup1982phase} and Gerchberg-Saxton \cite{gerchberg1972practical} methods remain to be state-of-the-art.  However, the re-weighted Amplitude flow method \cite{wang2017solving} appears to be highly competitive, if not practically identical in performance. 

 Convex methods appear to under-perform compared to their non-convex alternatives.
The convex method PhaseMax \cite{bahmani2017phase,goldstein2017icml,goldstein2016phasemax} seems to outperform the lifted convex relaxation PhaseLift (and its approximate low-rank version SketchyCGM \cite{yurtsever2017sketch}) on real image data. However, PhaseLift outperforms PhaseMax by a small margin on random Gaussian measurments.

\section{Discussion}

PhasePack was created with the goal of providing a unified framework in which researchers can investigate phase retrieval. By providing a common interface for different methods, and a simple API for testing methods on empirical datasets, we hope that PhasePack can help the community to better understand the strengths and weaknesses of different methods. 


\section*{Acknowledgments}

We would like to thank Ashok Veeraghavan and Chris Metzler at Rice University for making their empirical datasets available in PhasePack.  


\balance

\bibliographystyle{plain}
\bibliography{cites}

\end{document}